\documentclass[nopreprintline,12pt]{elsarticle}

 \usepackage{amsthm}
\usepackage{lineno}
\usepackage{latexsym,amscd,amsfonts,enumerate,supertabular}
\usepackage{tabularx}
\usepackage{array}
\usepackage{bigstrut}
\usepackage{graphicx}
\usepackage{epsfig}
\usepackage{amssymb}
\usepackage[intlimits]{amsmath}
\usepackage{algorithm}
\usepackage{xcolor}

\usepackage{multirow}
\usepackage{hhline}	
\usepackage{makecell}
\usepackage{color}
\usepackage{xcolor}
\usepackage{float}
\usepackage{graphicx,subfigure}
\usepackage{mathabx}

\usepackage{mathrsfs}

\definecolor{darkmagenta}{rgb}{0.55, 0.0, 0.55}

\newcolumntype{R}[1]{>{\raggedleft\let\newline\\\arraybackslash\hspace{0pt}}m{#1}}
\newcolumntype{L}[1]{>{\raggedright\let\newline\\\arraybackslash\hspace{0pt}}m{#1}}
\newcolumntype{C}[1]{>{\centering\let\newline\\\arraybackslash\hspace{0pt}}m{#1}}

\numberwithin{equation}{section}

\numberwithin{lemma}{section}
\newtheorem{theorem}[]{Theorem}
\numberwithin{theorem}{section}
\numberwithin{algorithm}{section}
\newtheorem{remark}[]{Remark}
\numberwithin{remark}{section}
\newtheorem{example}[]{Example}
\numberwithin{example}{section}
\newtheorem{definition}[]{Definition}
\numberwithin{definition}{section}

\numberwithin{corollary}{section}

\newcommand{\R}{\mathbb R}
\newcommand{\XX}{\mathfrak X}

\newcommand{\co}{\overline{\text{co}}}
\newcommand{\cl}{\text{cl}}

\begin{document}
\begin{frontmatter}
\title{{\bf
On strong law of large numbers for weakly stationary $\varphi$-mixing set-valued random variable sequences
}}
\author[ioit]{Luc T. Tuyen\corref{cor1}}
 \ead{tuyenlt@ioit.ac.vn}

\cortext[cor1]{Corresponding author}
\address[ioit]{Department of Data science and Big data analysis, Institute of Information Technology, Vietnam Academy of Science and Technology, Hanoi, Vietnam}
\begin{abstract}
\small  

In this paper we extend the notion of $\varphi$-mixing to set-valued random sequences that take values in the family of closed subsets of a Banach space.  
Several strong laws of large numbers for such $\varphi$-mixing sequences are stated and proved.  
Illustrative examples show that the hypotheses of the theorems are both natural and sharp.

\small  
\end{abstract}

\end{frontmatter}


\section{Introduction}\label{intro}
The law of large numbers plays a central role in probability theory and mathematical statistics~\cite{ibe2014fundamentals}.  
It has numerous practical applications, for example in finance~\cite{hull2018options}, economics~\cite{varian2014intermediate}, data mining and analysis~\cite{hastie2017elements}, and logistics~\cite{ballou2007business}.  
On a complete probability space $(\Omega,\mathcal{F},P)$, let $(x_1,x_2,\dots)$ be an independent and identically distributed (i.i.d.) sequence of real random variables with mean $\mu$ and finite variance $\sigma^{2}$.  
The classical law of large numbers asserts that $n^{-1}\sum_{k=1}^{n}x_k\to\mu$ as $n\to\infty$~\cite{ibe2014fundamentals}.  
Many extensions exist for non-identically distributed but independent sequences~\cite{geiringer1940generalization}.  

For \emph{set-valued} random sequences (random closed, typically convex sets) the law of large numbers in the i.i.d.\ case was established, for instance, in~\cite{artstein1975strong,van2012strong}.  
For arbitrary set-valued sequences that admit martingale-difference selections,~\cite{tuyen2020strong} obtained a corresponding result.  
However, the geometry of such set-valued sequences can still be highly dispersed even when they “oscillate’’ around $\{0\}$.

In practice, dependence is ubiquitous: Markov chains~\cite{iosifescu2014finite}, linear processes~\cite{characiejus2016weak}, and various mixing conditions~\cite{zhengyan1997limit}.  
Most of these notions of dependence admit a law of large numbers~\cite{kulik2017markov,zhengyan1997limit}.  
Specifically, for a sequence $\{x_n\}_{n\ge1}$, write $\mathcal F_m^{n}=\sigma(x_i:m\le i\le n)$ and define
$$
\varphi(n)=\sup_{m\ge1}\;
\sup_{A\in\mathcal F_1^{m},\,B\in\mathcal F_{m+n}^{\infty}}
|P(B\mid A)-P(B)|.
$$
If $\varphi(n)\to 0$ as $n\to\infty$, the sequence is called \emph{$\varphi$-mixing}.  
The following strong law for $\varphi$-mixing sequences appears in~\cite{gan2007strong,iosifescu1969random}.

\begin{theorem}\label{theo:1.1}\cite{gan2007strong}
	Let $\{x_n\}_{n\ge1}$ be a $\varphi$-mixing sequence with finite second moments and $\sum_{n=1}^{\infty}\varphi^{1/2}(n)<\infty$.  
	If $\{a_n\}_{n\ge1}$ is a non-decreasing sequence of positive numbers with $a_n\to\infty$ and
	$\sum_{n=1}^{\infty}\!Var(x_n)/a_n^{2}<\infty$, then
	\[
	\frac{\sum_{i=1}^{n}x_i-\sum_{i=1}^{n}E[x_i]}{a_n}\;\longrightarrow\;0
	\quad\text{a.s.}
	\]
\end{theorem}

Research on laws of large numbers for \emph{dependent} set-valued variables is still limited.  
In~\cite{fu2010note} a version for identically distributed $\varphi$-mixing random sets is stated, but the non-identical case remains open.  
Challenges include the fact that if a set-valued random variable has expectation $\{0\}$ it may degenerate to a single point~\cite{t2024representation,tuyen2022weak}, and that selections of non-identically distributed random sets may fail to preserve the dependence properties needed for single-valued laws of large numbers.

In this paper we introduce a new definition of $\varphi$-mixing for set-valued random sequences and prove several strong laws of large numbers under this dependence.

The remainder of the paper is organised as follows.  
Section~\ref{preliminaries} reviews notation and background on random sets.  
Section~\ref{result} states and proves our main results for weakly $\varphi$-mixing set-valued sequences, together with illustrative examples.  
Section~\ref{conc} concludes and outlines possible extensions.

\section{Notations and Preliminaries}\label{preliminaries}
Throughout this paper we work on a probability space
$\mathcal{P}=(\Omega,\mathcal{F},P)$ that satisfies the usual
conditions.  Let $(\mathfrak{X},\|\cdot\|_{\XX})$ be a separable Banach space,
and denote its dual by $\mathfrak{X}^{*}$ with the norm $\|\cdot\|_{\XX^*}$; the unit ball in
$\mathfrak{X}^{*}$ is written $S^{*}$.

Let $\mathbf K(\mathfrak X)$ be the family of non-empty closed subsets
of $\mathfrak X$.  Adding the subscripts $c$, $bc$, and $kc$ we obtain
the collections of non-empty closed \emph{convex}, closed \emph{bounded}
convex, and \emph{compact} convex subsets, respectively.
For $C\in\mathbf K(\mathfrak X)$ set
$\|C\|:=\sup\{\|x\|_{\XX}:x\in C\}$.
For $1\le p<\infty$ we write
$L^{p}[\Omega;\mathfrak X]$ for the Bochner space of measurable maps
$f:\Omega\to\mathfrak X$ with
\[
\|f\|_{p}:=\Bigl(\int_{\Omega}\|f(\omega)\|^{p}\,dP\Bigr)^{1/p}<\infty .
\]
(When $\mathfrak X=\mathbb R$ we simply write $L^{p}$.)

A set-valued map $F: \Omega \to \mathbf{K}(\mathfrak{X})$ is called a set-valued random variable or random set if for every closed subset $C$ of $\mathfrak{X}$, the set $\{\omega \in \Omega : F(\omega) \cap C \neq \emptyset\} \in \mathcal{F}$. A measurable function $f:\Omega\to\mathfrak X$ is a selection of $F$ when $f(\omega)\in F(\omega)$ for all $\omega\in\Omega$. The function $f$ is termed an \emph{almost everywhere selection} of $F$ if $f(\omega) \in F(\omega)$ for almost every $\omega \in \Omega$. We also denote the collection of all set-valued random variables by $\mathcal{U}[\Omega, \mathcal{F}, P; \mathbf{K}(\mathfrak{X})]$.\\
For every number $1 \leq p \leq \infty$, we denote
$$ S_F^p(\mathcal{F}) = \left\{ f \in L^p[\Omega; \mathfrak{X}] : f(\omega) \in F(\omega), \text{ a.s.} \right\} $$
as the set of all $p$-order integrable selections of the set-valued random variable $F$. For simplicity, we write $S_F$ instead of $S^1_F$. Note that $S_F^p$ is a closed subset of $L^p[\Omega; \mathfrak{X}]$.

A set-valued random variable $F$ is said to be integrable if $S^1_F$ is non-empty. $F$ is called $L^p$-integrably bounded (or strongly integrable) if there exists $\rho \in L^p[\Omega, \mathcal{F}, \mathbb{R}]$ such that $|x| \le \rho(\omega)$ for all $x \in F(\omega)$ and for all $\omega \in \Omega$. We denote $L^p[\Omega, \mathcal{F}, P; \mathbf{K}(\mathfrak{X})]$ as the set of all $L^p$-integrably bounded set-valued random variables.

For a set-valued random variable $F$, the Aumann integral of $F$ is defined by 
$$ E[F] = \int_{\Omega}FdP = \left\{\int_{\Omega} fdP : f \in S_F \right\}, $$
where $\int_{\Omega} fdP$ is the conventional Bochner integral in $L^1[\Omega; \mathfrak{X}]$. Since $\int_{\Omega} FdP$ is generally not a closed set, except under certain conditions such as when $\mathfrak{X}$ has the Radon-Nikodym property and $F \in L^1[\Omega; K_{kc}(\mathfrak{X})]$, or when $\mathfrak{X}$ is reflexive and $F \in L^1[\Omega; K_{c}(\mathfrak{X})]$ (note that $L^p[\Omega; \mathfrak{X}]$ for $1 < p < \infty$ and all finite-dimensional spaces are reflexive).

For a non–empty closed convex set $C\subset \mathfrak{X}$ we write
$$
V_{\infty}(C):=\bigl\{u\in \mathfrak{X} \;:\; C+\lambda u\subseteq C\;\text{ for every } \lambda\ge 0\bigr\}
$$
and call $V_{\infty}(C)$ the \emph{recession cone} of $C$ (see \cite[page 50]{bertsekas2003convex}).

We denote $\mathcal{P}_0(\mathfrak{X})$ as the collection of non-empty subsets of $\mathfrak{X}$. For all $A, B \in \mathcal{P}_0(\mathfrak{X})$ and $\lambda \in \mathbb{R}$, we define

$$\begin{aligned}
	A+B=\{a+b:a\in A, b\in B\}\\
	\lambda A=\{\lambda a:a\in A\}.
\end{aligned}$$

Note that if $A, B \in \mathbf{K}_k(\mathfrak{X})$, then $A + B \in \mathbf{K}_k(\mathfrak{X})$.

For every $A\in\mathcal{P}_0(\mathfrak X)$ and $x \in \mathfrak{X}$, the distance between $x$ and $A$ is defined by $d(x,A)=\inf_{y\in A}d(x,y)$, where $d(x,y)=\|x-y\|_{\mathfrak X}$.

The Hausdorff distance on $\mathcal{P}_0(\mathfrak X)$ is defined as follows:
$$H(A,B)=\max\{\sup_{a\in A}d(a,B),\sup_{b\in B}d(b,A)\}.$$
The Hausdorff distance between $A$ and $\{0\}$ is denoted by $\|A\|=H(A,0)$. For a finite set $E$, we write $|E|$ for its cardinality, $|\cdot|$ is otherwise used only for absolute values of scalars.

For each $A\in\mathbf{K}(\mathfrak X)$, the support function of $A$ is defined by $s(x^*,A)=\sup_{a\in A}\langle x^*,a \rangle$ for each $x^*\in\mathfrak{X}^*$. Note that $s(x^*,\cl(A+B))=s(x^*,A+B)=s(x^*,A)+s(x^*,B)$ and $s(x^*,\lambda A)=\lambda s(x^*,A)$ for all $x^*\in\mathfrak X^*$ and $\lambda \ge 0$.

Additionally, we denote $\co{A}$ as the closed convex hull of $A$. Then, $x \in \co{A}$ if and only if $\langle x^*,x\rangle \le s(x^*,A)$ for all $x^* \in \mathfrak{X}^*$.

Let $\{A_n, A\}$ be a sequence of closed subsets of $\mathfrak{X}$. $A_n$ is said to converge to $A$ in the Hausdorff sense if $\lim\limits_{n\to\infty}H(A_n,A)=0$. $A_n$ is said to converge to $A$ in the Kuratowski-Mosco sense \cite{mosco1969convergence} if 
$$\mathrm{w}\text{-}\limsup_{n\to\infty}{A_n}=A=\mathrm{s}\text{-}\liminf_{n\to\infty}{A_n}.$$

Here, $$\mathrm{w}\text{-}\limsup_{n\to\infty}{A_n}=\{x=\text{w-}\lim_{k\to\infty}a_{k}:a_{k}\in A_{n_k}, \{A_{n_k}\}\subset \{A_n\},n,k\ge 1\},$$
 and $$\mathrm{s}\text{-}\liminf_{n\to\infty}{A_n}=\{x=\text{s-}\lim_{n\to\infty}a_{n}:a_{n}\in A_{n}, n\ge 1\}.$$
 
 Furthermore, $\text{s-}\lim\limits_{n\to\infty}x_n=x$ means that $\|x_n-x\|_{\XX}\to 0$ and $\text{w-}\lim\limits_{n\to\infty}x_n=x$ means that $x_n$ converges weakly to $x$.

\section{Main results}\label{result}
In this section we introduce the notions of \emph{weak stationarity} and
$\varphi$-mixing for sequences of set-valued random variables, and we
prove several strong laws of large numbers that describe their
convergence both in the Hausdorff metric and in the
Kuratowski–Mosco sense.

Let $\{X_n, n \ge 1\}$ be a sequence of set-valued random variables defined on $\mathcal{U}\left[ \Omega ,\mathcal{F},P;\mathbf{K}\left( \mathfrak{X} \right) \right]$. Similar to the Introduction section, let $\mathcal{F}_{m}^{n}=\sigma \left( {{X}_{i}},m\le i\le n \right)$ for all $m,n\in \mathbb{N}$, which is the $\sigma$-algebra generated by $X_m,X_{m+1},...,X_n$. For two $\sigma$-algebras $\mathcal{A}, \mathcal{B} \in \mathcal{F}$, define
$$\varphi \left( \mathcal{A},\mathcal{B} \right)=\underset{ A\in \mathcal{A},B\in \mathcal{B},P(A)\ne 0}{\mathop{\sup }}\,\left| P(B|A)-P(B) \right|.$$

We define the dependence coefficient $\varphi$ as follows:
$$\varphi (n)=\underset{k\in \mathbb{N}}{\mathop{\sup }}\,\varphi \left( \mathcal{F}_{1}^{k},\mathcal{F}_{k+n}^{\infty } \right),\,\,n\ge 0.$$

\begin{definition}
	A sequence of set-valued random variables $\{X_n, n \ge 1\}$ is called a $\varphi$-mixing sequence of set-valued random variables if $\varphi(n) \to 0$ as $n \to \infty$.
\end{definition}

\begin{definition}\label{def3.2}
	A sequence of set-valued random variables $\{X_n, n \ge 1\}$ is called weakly stationary if $E(X_n) = A$ for all $n \in \mathbb{N}$, where $A$ is a non-empty closed subset of $\mathfrak{X}$.
\end{definition}

Note that, regarding the definition of a strictly stationary sequence of set-valued random variables proposed by Wang \cite{wang1997set}, this definition is weaker. According to Wang, a sequence of set-valued random variables $\{X_1, X_2, \ldots\}$ is called strictly stationary if for every $(\mathcal{U}_1, \mathcal{U}_2, \ldots, \mathcal{U}_i) \subset \mathcal{F}$ and $(t_1, t_2, \ldots, t_i) \subset \mathbb{N}$ and $t\in \mathbb{N}$, then
$$\begin{aligned}
	P\{\omega:X_{t_k}(\omega)\in \mathcal U_{k}, 1\le k\le i\}=\\
	P\{\omega:X_{t_k+t}(\omega)\in \mathcal U_{k}, 1\le k\le i\}
\end{aligned} .$$

Thus, the concept of a strictly stationary set-valued random variable is expressed in terms of translation invariant probability distribution. Because strict stationarity implies identical distributions, the sequence automatically satisfies $E[X_n]=E[X_1]$ for all $n\in\mathbb N$. Therefore, the properties of weakly stationary sequences of set-valued random variables will also hold for strictly stationary sequences. We can provide an example demonstrating a weakly stationary sequence of set-valued random variables that is not strictly stationary.

\begin{example}
	In the Banach space $\mathbb{R}^2$, consider the sequence of set-valued random variables $X_n = \{(x_1, x_2) : x_1^2 + x_2^2 \le r_n^2\}$, which represents a random ball centered at $O$ with radius $r_n$. Here, $r_n$ is a random variable that depends on $n$ defined as follows:
	\begin{itemize}
		\item For $n$ is an even number, $r_n$ is uniformly distributed on $[0.9, 1.1]$.
		\item For $n$ is an odd number, $r_n$ follows a normal distribution $N(1, 0.1)$.
	\end{itemize}
	It is clear that $\{X_n\}_{n \ge 1}$ is weakly stationary because $E[X_n]$ is a ball centered at $O$ with radius $1$ for all $n \ge 1$, but it is not strictly stationary due to the differing distributions between even and odd $n$.
\end{example}

Next, the author provides an example of a weakly stationary $\varphi$-mixing sequence of set-valued random variables that takes values on compact convex subsets of $\mathbb{R}$ without degenerating into single-valued variables.

\begin{example}\label{ex:3.2}
	Suppose $\{x_n\}_{n\ge 1}$ is a bounded $\varphi$-mixing sequence of random variables with a common expectation of $\mu$ and variance of $\sigma^2$. We define
	$$ X_n=[x_n,x_n+1].$$
	This is a random line segment in $\mathbb{R}$, and each $X_n$ is a compact convex set. Consequently, $\{X_n\}_{n \ge 1}$ is a weakly stationary and $\varphi$-mixing sequence of set-valued random variables.
\end{example}

Indeed, it is clear that 
$E[X_n]=[E[x_n],E[x_n+1]]=[\mu,\mu+1]=A$, which is constant with respect to $n$.\\
On the other hand, since $\{x_n,n \ge 1\}$ is a $\varphi$-mixing sequence of random variables, $\{X_n\}_{n \ge 1}$ is also a $\varphi$-mixing sequence of set-valued random variables.

\begin{remark}\label{re:3.1}
	If the sequence $\{X_n;n \ge 1\}$ is a weakly stationary sequence of set-valued random variables, then the sequence $\{s(x^*,X_n)\}_{n\ge 1}$ is also a weakly stationary sequence of single-valued random variables for all $x^* \in \mathfrak{X}^*$ due to $E[s(x^*, X_n)] = s(x^*, E[X_n]) = s(x^*, A)$. Furthermore, if $\{X_n\}_{n \ge 1}$ is a $\varphi$-mixing sequence of set-valued random variables, then $\{s(x^*, X_n)\}_{n \ge 1}$ is also a $\varphi$-mixing sequence of random variables for all $x^* \in \mathfrak{X}^*$ because the $\sigma$-algebra generated by $s(x^*, X_n)$ is a sub-$\sigma$-algebra of the $\sigma$-algebra generated by $X_n$.
\end{remark}

\begin{theorem}\label{theo:3.1}
	Let $\{X_n, n \ge 1\}$ be a weakly stationary and $\varphi$-mixing sequence in ${{L}^{2}}\left[ \Omega,\mathcal{F},P;{{\mathbf{K}}_{kc}}\left( \mathfrak{X} \right) \right]$ with $E[X_n] = A$ for all $n \ge 1$, and the following conditions are satisfied:\\
	(i) $\sum\limits_{n=1}^{\infty}{{\varphi }^{1/2}}(n)<\infty $,\\
	(ii) $\sum\limits_{n=1}^{\infty}\frac{E|s(x^*,{X}_{n})-s(x^*, A)|^2}{n^2}<\infty $ for all $x^*\in S^*$.\\
	Then
	$$H\left( \frac{1}{n}\sum\limits_{k=1}^{n}{X_k};A \right)=0,\,a.s.$$	
\end{theorem}

\begin{proof}
	According to Corollary 1.1.10 (from \cite{li2013limit}), for compact convex random variables, we have 
	
	\begin{equation}\label{eq:3.1}
		H\left( \frac{1}{n}\underset{k=1}{\overset{n}{\mathop \sum }}\,{{X}_{k}},A \right)=\left\|\frac{1}{n}\underset{k=1}{\overset{n}{\mathop \sum }}\,s(\cdot ,{{X}_{k}})-s(\cdot ,A)\right\|_{C(S^*,d_{w}^{*})},
	\end{equation}	
	where $C(S^*,d_{w}^{*})$ is the space of bounded functions in $S^*$ with the weak metric $d_w^*$ defined by $d_w^*(x_1^*,x_2^*)=\sum_{i=1}^{\infty}\frac1{2^i}|\langle x_1^*,x_i\rangle-\langle x_2^*,x_i\rangle|$ and $x_1, x_2, \ldots$ being dense in the closed unit ball of $\mathfrak{X}$.
	
	Furthermore, for every $x^*\in {{S}^{*}}$, $s(x^*,{{X}_{1}}),\cdots ,s(x^*,{{X}_{n}})$ is a weakly stationary $\varphi$-mixing sequence with the common expectation at $s(x^*, A)$, and
	$ E[s(x^*,{{X}_{n}})]=s(x^*,E[{{X}_{n}}])$ $=s(x^*,clE[{{X}_{n}}])=s(x^*,A) $ 
	for all $x^* \in S^*$.
	
	Therefore, by the law of large numbers for $\varphi$-mixing random sequences (Theorem \ref{theo:1.1}), we have
	$$ \left| \frac{1}{n}\underset{k=1}{\overset{n}{\mathop \sum }}\,s({{x}^{*}},{{X}_{k}})-s({{x}^{*}},A) \right|\to 0 $$ 
	as $n$ approaches infinity for all $x^* \in S^*$.
	
	This demonstrates that the right-hand side of \eqref{eq:3.1} converges to $0$ almost surely as $n$ goes to infinity.
\end{proof}

\begin{example}\label{ex:3.3}
	Example \ref{ex:3.2} with the assumption that $\{x_n,n\ge 1\}$ have dependence coefficients strong enough to ensure that $\sum\limits_{n=1}^{\infty}\varphi^{\frac12}(n)<\infty$, is an example of a sequence of set-valued random variables that satisfies the conditions of Theorem \ref{theo:3.1}. Furthermore, we can directly calculate the limit $\lim\limits_{n\to\infty}H\left( \frac{1}{n}\sum\limits_{k=1}^{n}{{{X}_{k}}};A \right)=0,\,a.s.$ without applying the results of Theorem \ref{theo:3.1}.
\end{example}

Indeed, condition (i) of Theorem \ref{theo:3.1} is satisfied by the assumptions. \\
For $x^* \ge 0$, we have $s(x^*,X_n)=x^*(x_n+1)$ and , $s(x^*,A)=x^*(\mu+1)$ (where $A = [\mu, \mu + 1]$). Thus,
$$E|s(x^*,X_n)-s(x^*,A)|^2=E|x^*(x_n+1)-x^*(\mu+1)|^2=(x^*)^2E[x_n-\mu]^2=\sigma^2.$$
For $x^* < 0$, we have $s(x^*,X_n)=x^*(x_n)$ and $s(x^*,A)=x^*\mu$. Therefore,
$$E|s(x^*,X_n)-s(x^*,A)|^2=E|x^*(x_n-\mu)|^2=(x^*)^2E[x_n-\mu]^2=\sigma^2.$$
Thus,
$$\sum_{n=1}^{\infty}\frac{E|s(x^*,X_n)-s(x^*,A)|^2}{n^2}=\sigma^2\sum_{n=1}^{\infty}\frac{1}{n^2}<\infty.$$
Therefore, condition (ii) of Theorem \ref{theo:3.1} is satisfied.

We have:
$$\frac1n\sum_{k=1}^nX_k=\left[\frac1n\sum_{k=1}^nx_n;\frac1n\sum_{k=1}^n(x_n+1)\right].$$
According to the strong law of large numbers for the $\varphi$-mixing sequence $\{x_n, n \ge 1\}$ (Theorem \ref{theo:1.1}), $\frac1n\sum_{k=1}^nx_n\to \mu \;a.s.$ and $\frac1n\sum_{k=1}^n(x_n+1)=\frac1n\sum_{k=1}^nx_n+1\to \mu+1 \;a.s.$ as $n\to\infty$.\\
Thus, for all $a\in \left[\frac1n\sum_{k=1}^nx_n;\frac1n\sum_{k=1}^n(x_n+1)\right]$, we have $\inf\limits_{b\in [\mu,\mu+1]}d(a,b)\to 0$ as $n\to \infty$.\\
Similarly, for every $b \in [\mu, \mu + 1]$ and every $\varepsilon > 0$, there exists a sufficiently large number $N > 0$ such that $d(b, a) < \varepsilon$.\\
Therefore, $\inf\limits_{a\in \left[\frac1n\sum_{k=1}^nx_n;\frac1n\sum_{k=1}^n(x_n+1)\right]}d(b,a)\to 0$ as $n\to \infty$.\\
Thus, $\lim_{n\to\infty}H\left(\frac1n\sum_{k=1}^nX_n,A\right)=0,\; a.s.$

Now we will extend the strong law of large numbers for weakly stationary $\varphi$-mixing sequences of set-valued random variables taking values in the space of closed compact  (not necessarily convex) sets.

\begin{theorem}\label{theo:3.2}
	Let $\{X_n,n\ge 1\}$ be a weakly stationary $\varphi$-mixing set-valued random variable sequence in $L^2[\Omega,\mathcal F,P;\mathbf{K}_k(\mathfrak{X})]$ with $E[X_n]=A$ for all $n\ge 1$ and the following conditions are satisfied:\\
	(i) $\sum\limits_{n=1}^{\infty}{{\varphi }^{1/2}}(n)<\infty $ and\\
	(ii) $\sum\limits_{n=1}^{\infty}\frac{E|s(x^*,\co{X}_{n})-s(x^*,\co A)|^2}{n^2}<\infty $ for all $x^*\in S^*$\\
	Then
	$$\lim_{n\to\infty}H\left( \frac{1}{n}\sum\limits_{k=1}^{n}{{X_k}};\co A \right)=0,\,a.s.$$
\end{theorem}

\begin{proof}
	By Theorem \ref{theo:3.1}, we have
	\begin{equation}\label{eq:3.2}
		\lim_{n\to\infty}H\left(\frac1n\sum_{k=1}^n\co{X_k},\co{A}\right)=0\, a.s.
	\end{equation}
	
	By Lemma 3.1.4 in \cite{li2013limit}, the proof of Theorem \ref{theo:3.2} is completed.
\end{proof}

\begin{example}\label{ex:3.4}
	Let $\{x_n, n \geq 1\}$ be a sequence of $\varphi$-mixing random variables that are identically distributed, with mean $0$ and variance $1$, satisfying $\sum\limits_{n=1}^\infty \varphi^{\frac{1}{2}}(n) < \infty$. Define a sequence of random sets $\{X_n, n \geq 1\}$ taking values in $\mathcal{P}(\mathbb{R})$ as follows: 
	$$X_n = \{x_n, x_n + 1\}.$$
	Then $\{X_n, n \geq 1\}$ satisfies the conditions in Theorem \ref{theo:3.2}, and the strong law of large numbers can be proved directly.
\end{example}

Indeed, since the sigma-algebra generated by $\{X_i, m \leq i \leq n\}$ is equivalent to the sigma-algebra generated by $\{x_i, m \leq i \leq n\}$ for all $n \geq m \geq 1$, it is clear that $\{X_n, n \geq 1\}$ is compact, non-convex, and weakly stationary, with $A= E[X_n] = \{0, 1\}$ and $\varphi$-mixing satisfying condition (i).

Condition (ii), which states that 
$$\sum\limits_{n=1}^\infty \frac{E|s(x^*, \co{X_n}) - s(x^*, \co{A})|^2}{n^2} < \infty \quad \text{for all } x^* \in S^*,$$ 
is verified similarly to Example \ref{ex:3.3}.

We have 
$$\frac{1}{n}\sum\limits_{k=1}^n X_k = \left\{\frac{1}{n}\sum\limits_{k=1}^n x_k + \frac{i}{n} : 0 \leq i \leq n \right\},$$ 
which includes $n+1$ points equally spaced by $\frac{1}{n}$ from $\frac{1}{n}\sum\limits_{k=1}^n x_k$ to $\frac{1}{n}\sum\limits_{k=1}^n x_k + 1$. Here, $\co{A} = \co{\{0, 1\}} = [0, 1]$. 

It follows that 
$$\sup\limits_{a \in \frac{1}{n}\sum\limits_{k=1}^n X_k} d(a, [0, 1]) = \frac{1}{n}\sum\limits_{k=1}^n x_k$$ 
and 
$$\sup\limits_{b \in [0, 1]} d(b, \frac{1}{n}\sum\limits_{k=1}^n X_k) = \frac{1}{n}\sum\limits_{k=1}^n x_k + \frac{1}{2n}.$$ 

Thus, 
$$H\left(\frac{1}{n}\sum\limits_{k=1}^n X_k, [0, 1]\right) = \frac{1}{n}\sum\limits_{k=1}^n x_k + \frac{1}{2n}.$$ 

By the strong law of large numbers for $\{x_n : n \geq 1\}$, $\frac{1}{n}\sum\limits_{k=1}^n x_k \to 0$ as $n \to \infty$, and since $\frac{1}{2n} \to 0$, we conclude that 
$$H\left(\frac{1}{n}\sum\limits_{k=1}^n X_k, [0, 1]\right) \to 0 \quad \text{a.s.}$$

Now, we proceed to establish the strong law of large numbers for the weakly stationary $\varphi$-mixing sequence of set-valued random variables in the sense of Kuratowski-Mosco.

\begin{theorem}\label{theo:3.3}
	Let $\{X_n, n\ge1\}$ be a weakly stationary $\varphi$-mixing sequence of
	set–valued random variables in
	$\mathcal U[\Omega,\mathcal F,P;\mathbf K(\mathfrak X)]$ with
	$E[X_n]=A$ for all $n\ge1$.
	Assume
	
	\begin{itemize}
		\item[(i)] $\displaystyle\sum_{n=1}^{\infty}\varphi^{1/2}(n)<\infty$;
		\item[(ii)] For every $a\in A$, there exists a square-integrable selection $\{x_n\}\subset S_{X_n}(\mathcal F)$ satisfying $E[x_n]=a$ for all $n\ge 1$ and
		$$
		\sum_{n=1}^{\infty}
		\frac{E\|x_n-a\|^{2}_{\XX}}{n^{2}}<\infty;
		$$
		\item [(iii)] For every $x^*\in S^*$ with $s(x^*,A)<\infty$,
		\[
		\sum_{n=1}^{\infty}
		\frac{E\,\big|\,s(x^*,X_n)-s(x^*,A)\,\big|^{2}}{n^{2}} \;<\;\infty .
		\]
		
	\end{itemize}
	
	Then
	$$
	S_n \xrightarrow{\,\mathrm{K\text{-}M}\,} \;D.
	$$
	{Here} $S_n = \frac{1}{n} \mathrm{cl} \sum_{i=1}^n X_i$ and $D=\overline{\operatorname{co}}A.$
\end{theorem}

\begin{proof}
	
	First, we prove $D\subset\displaystyle\mathrm{s}\text{-}\liminf_{n\to\infty}S_n$. 
	Recall that $D = \co{A}$ and $S_n = \frac{1}{n} \cl \sum_{i=1}^n X_i$. 
	Fix $d\in D=\overline{\operatorname{co}}A$ and $\varepsilon>0$.
	Choose $a_1,\dots,a_r\in A$ and rational weights
	$\lambda_j=p_j/q$ with $p_j\in\mathbb N$, $\sum_{j=1}^r p_j=q$, such that
	\begin{equation}\label{eq:d_approx}
	\Bigl\|\sum_{j=1}^{r}\lambda_j a_j - d\Bigr\|_{\XX}<\varepsilon .
	\end{equation}
	Let $L_1,\dots,L_q\in\{1,\dots,r\}$ be a $q$–periodic label sequence with
	$\#\{t\in\{1,\dots,q\}:L_t=j\}=p_j$ for every $j$.
	
	For each $j=1,\dots,r$, assumption (ii) (applied to $a=a_j$) yields a
	global square–integrable selection $\{x^{(j)}_n\}_{n\ge1}\subset S^{2}_{X_n}(\mathcal F)$ such that
	\[
	E[x^{(j)}_n]=a_j\quad\forall n\ge1,
	\qquad
	\sum_{n=1}^{\infty}\frac{E\|x^{(j)}_n-a_j\|^2_{\XX}}{n^2}<\infty .
	\]
	Since $\sigma(x^{(j)}_n)\subseteq\sigma(X_n)$, each sequence $\{x^{(j)}_n\}$ inherits
	the $\varphi$–mixing property with coefficients bounded by those of $\{X_n\}$.
	
	\medskip
	Define a single selection $\{x_n\}_{n\ge1}$ by the $q$–periodic schedule
	\[
	\quad x_n:=x^{(L_{((n-1)\bmod q)+1})}_n,\qquad n=1,2,\dots
	\]
	Then $x_n\in X_n$ for every $n$ (the index on the chosen selection matches $n$),
	and $\{x_n\}$ is again $\varphi$–mixing.
	
	\medskip
	Write $N=qk(N)+\ell(N)$ (we abbreviate $k:=k(N)$ and $\ell:=\ell(N)$ to write $N=qk+\ell$) with $k\ge0$ and $0\le\ell<q$. For each $j=1,\dots,r$ set
	$I_j(N):=\{1\le n\le N:\ L_{((n-1)\bmod q)+1}=j\}$; then
	$|I_j(N)|=k\,p_j+r_j$ with $0\le r_j\le p_j$. Since $k\to\infty$ as $N\to\infty$, $\frac{|I_j(N)|}{N}=\frac{kp_j+r_j}{kq+\ell}=\frac{p_j}q.\frac{k}{k+\ell/q}+\frac{r_j}{N}\to p_j/q=\lambda_j$ as $N\to\infty$.
		
	For each $j=1,2,...,r$ and each $t\in\{1,\dots,q\}$ with $L_t=j$,
	define the arithmetic–progression subsequence
	\[
	\xi^{(j,t)}_b:=x^{(j)}_{(b-1)q+t},\qquad b=1,2,\dots
	\]
	Then $\{\xi^{(j,t)}_b\}_{b\ge1}$ is $\varphi$–mixing with coefficients
	{$\varphi(q(b-1))$}, hence
	{$\sum_{b\ge1}\varphi(q(b-1))^{1/2}\le \sum_{n\ge1}\varphi(n)^{1/2}<\infty$}
	(since $\varphi$ is nonincreasing). Moreover,
	\begin{equation}\label{eq:xi_b_bounded}
	\sum_{b=1}^{\infty}\frac{E\|\xi^{(j,t)}_b-a_j\|^2_{\XX}}{b^2}
	= \sum_{b=1}^{\infty}\frac{E\|x^{(j)}_{(b-1)q+t}-a_j\|^2_{\XX}}{b^2}
	\;\le\; q^{2}\sum_{n=1}^{\infty}\frac{E\|x^{(j)}_{n}-a_j\|^2_{\XX}}{n^2}
	<\infty .
	\end{equation}
	Since $\mathfrak{X}$ is separable, take $\mathcal G=\{y_m^*\}_{m\ge1}\subset S^*$ be a fixed countable dense set.
	For each $m$ and each arithmetic progression $(j,t)$ used in the construction, set $\eta^{(j,t)}_b(y^{*}_m)
	:=\langle y^{*}_m,\,\xi^{(j,t)}_b-a_j\rangle$. Then
	$\sum_{b\ge1}\frac{E|\eta^{(j,t)}_b(y^{*}_m)|^{2}}{b^{2}}
	\le \sum_{b\ge1}\frac{E\|\xi^{(j,t)}_b-a_j\|^{2}_{\XX}}{b^{2}}<\infty$ by \eqref{eq:xi_b_bounded},
	and the subsequence $\{\eta^{(j,t)}_b(y^{*}_m)\}$ is also $\varphi$–mixing with
	$\sum_{b\ge1}\varphi(qb)^{1/2}\le\sum_{n\ge1}\varphi(n)^{1/2}<\infty$. 
	Theorem~\ref{theo:1.1} yields an event $\Omega_{m}^{(j,t)}$ with probability~1 on which
	\[
	\frac1k\sum_{b=1}^{k}\langle y_m^*,\,\xi^{(j,t)}_b-a_j\rangle=\frac1k\sum_{b=1}^{k}\eta^{(j,t)}_b(y^{*}_m) \longrightarrow 0 .
	\]
	Set $\displaystyle \Omega_0:=\bigcap_{(j,t)}\ \bigcap_{m=1}^{\infty} \Omega_{m}^{(j,t)}$.
	Then $P(\Omega_0)=1$.
	Hence, for every $\omega\in\Omega_0$ and every finite set $F\subset\mathcal G$,
	\begin{equation}\label{eq:3.4}
	\max_{y^*\in F}\left|
	\frac1k\sum_{b=1}^{k}\langle y^*,\,\xi^{(j,t)}_b-a_j\rangle
	\right|\ \longrightarrow\ 0 .
	\end{equation}
		
	To pass from scalar to norm convergence, fix $\delta\in(0,1)$. For each $k$, let
	$V_k:=\mathrm{span}\{\xi^{(j,t)}_1-a_j,\dots,\xi^{(j,t)}_k-a_j\}$ and select a finite $\delta$–net $F_k\subset\mathcal G$ such that
	$\{\,y^*|_{V_k}:y^*\in F_k\,\}$ is {a} $\delta$–net of the unit sphere of $V_k^*$,  {i.e., every element of the unit sphere has distance less than $\delta$ to some
	$y^* \in F_k$} (by Hahn–Banach extension and the density of $\mathcal G$ in $S^*$ this is always possible).
	The standard net estimate then yields, for every $v\in V_k$,
	\begin{equation}\label{eq:3.5}
	(1-\delta)\,\|v\|_{\XX}\ \le\ \max_{y^*\in F_k} \big|\langle y^*,v\rangle\big|.
	\end{equation}
	(The estimate follows by taking $x^*_k\in V_k^*$ with
	$\|x^*_k\|_{V_k^*}=1$ and $x^*_k(v)=\|v\|_{\mathfrak X}$ for $v\in V_k$,
	then choosing $y^*\in F_k$ with { $\|\,y^*|_{V_k}-x^*_k\,\|_{V_k^*}\le\delta\|v_k\|_{\mathfrak X}$ }.
	Here $\|\cdot\|_{V_k^*}$ denotes the operator norm on $V_k^*$, i.e.
	$\|f\|_{V_k^*}=\sup\{\,|f(w)|:\ w\in V_k,\ \|w\|_{\mathfrak X}\le1\,\}$).
		
	Because $F_k$ varies with $k$, we cannot invoke \eqref{eq:3.4} with $F=F_k$
	uniformly in $k$. We therefore \emph{freeze} the family by setting
	\[
	H_m:=\bigcup_{i=1}^{m} F_i \qquad(\text{finite for each }m).
	\]
	
	Fix $\omega\in\Omega_0$. By \eqref{eq:3.4}, for each $m$
	there exists $K_m=K_m(\omega)$ such that for all $k\ge K_m$,
	\[
	\max_{y^*\in H_m}\Bigl|\frac1k\sum_{b=1}^k\langle y^*,\xi_b^{(j,t)}-a_j\rangle\Bigr|
	\le \frac1m.
	\]
	Choose inductively $k_1\ge K_1$ and $k_r\ge\max\{k_{r-1}+1,\;K_r\}$ for $r\ge2$.
	Then $k_r\uparrow\infty$ and, at $k=k_r$,
	\[
	\max_{y^*\in H_{k_r}}
	\Bigl|\frac1{k_r}\sum_{b=1}^{k_r}\langle y^*,\xi_b^{(j,t)}-a_j\rangle\Bigr|
	\le \frac1{k_r}\xrightarrow[r\to\infty]{}0.
	\]
	
	As $F_{k_r}\subset H_{k_r}$, we obtain
	\[
	\max_{y^*\in F_{k_r}}\Bigl|\frac1{k_r}\sum_{b=1}^{k_r}\bigl\langle y^*,\,\xi^{(j,t)}_b-a_j\bigr\rangle\Bigr|
	\longrightarrow 0 .
	\]
	With $v_{k}:=\frac1{k}\sum_{b=1}^{k}(\xi^{(j,t)}_b-a_j)\in V_k$, inequality \eqref{eq:3.5}
	gives along the subsequence $\{k_r\}$:
	\[
	\|v_{k_r}\|_{\mathfrak X}
	\ \le\ \frac{1}{1-\delta}\,
	\max_{y^*\in F_{k_r}}\Bigl|\frac1{k_r}\sum_{b=1}^{k_r}\bigl\langle y^*,\,\xi^{(j,t)}_b-a_j\bigr\rangle\Bigr|
	\ \xrightarrow[r\to\infty]{}\ 0 .
	\]
	Because $\delta\in(0,1)$ is arbitrary, 
	\begin{equation}\label{eq:3.6prime}
		\frac1{k_r}\sum_{b=1}^{k_r}\xi^{(j,t)}_b \ \xrightarrow[r\to\infty]{a.s.}\ a_j .
	\end{equation}
	
	Choose an increasing sequence $N_r\uparrow\infty$ such that $k(N_r)=k_r$ (i.e. $N_r=qk_r$ and $\ell(N_r)=0$).
	Note that the set $\{t:L_t=j\}$ has cardinality $p_j$ independent of $N$.
	Hence, along $N_r$ with $k(N_r)=k_r$,
	\begin{equation}\label{eq:mid-term}
	\frac{1}{p_j}\sum_{t:\,L_t=j}\Bigl(\frac{1}{k_r}\sum_{b=1}^{k_r}\xi^{(j,t)}_b\Bigr)
	\ \xrightarrow[r\to\infty]{a.s.}\ 
	\frac{1}{p_j}\sum_{t:\,L_t=j} a_j \;=\; a_j,
	\end{equation}
	
		On the other hand, for $|I_j(N_r)|=k_rp_j+r_j$ then $I_j(N_r)=\bigcup_{t:\,L_t=j}\{(b-1)q+t:1\le b\le k_r\}\ \cup\ R_j(N_r)$
	with $R_j(N_r)\subset\{k_rq+1,\dots,k_rq+\ell\}$ and $|R_j(N_r)|=r_j$, we have
	\begin{equation}\label{eq:Ij-decomp}
		\frac{1}{|I_j(N_r)|}\sum_{n\in I_j(N_r)} x^{(j)}_n
		= \frac{k_r p_j}{k_r p_j + r_j}\cdot \frac{1}{p_j}\sum_{t:\,L_t=j}
		\Bigl(\frac{1}{k_r}\sum_{b=1}^{k_r}\xi^{(j,t)}_b\Bigr)
		\ +\ \frac{1}{|I_j(N_r)|}\sum_{n\in R_j(N_r)} x^{(j)}_n .
	\end{equation}
	
	By \eqref{eq:mid-term} and $\frac{k_r p_j}{k_r p_j + r_j}\to1$ as $r\to\infty$, the middle term of \eqref{eq:Ij-decomp} tends to $a_j$ a.s.. { By (ii), Borel–Cantelli and Chebyshev, } $\|x_n^{(j)}\|_{\XX}/n\to0$, hence for any $\varepsilon>0$ and $N_r$ large enough, $\|x_n^{(j)}\|_{\XX}\le \varepsilon.(k_rq+q)$ for all $n\in R_j(N_r)$, which yields 
	
	\[\frac{1}{|I_j(N_r)|}\sum_{n\in R_j(N_r)} \|x^{(j)}_n\|_{\XX}\le \tfrac{r_j\varepsilon (k_rq+q)}{k_rp_j}\le \tfrac{k_r+1}{k_r}\varepsilon q \xrightarrow[r\to\infty]{}\ 0.\]
	Consequently,
	\begin{equation}\label{eq:Ij-avg-subseq}
		\frac{1}{|I_j(N_r)|}\sum_{n\in I_j(N_r)} x^{(j)}_n\ \xrightarrow[r\to\infty]{\ \|\cdot\|_{\mathfrak X}\ }\ a_j
		\quad \text{a.s.}
	\end{equation}
	
	Since also $|I_j(N_r)|/N_r\to p_j/q=\lambda_j$, we obtain
	\begin{equation}\label{eq:blocksum-subseq}
		\frac{1}{N_r}\sum_{n=1}^{N_r} x_n
		=\sum_{j=1}^{r}\frac{|I_j(N_r)|}{N_r}\,
		\Big(\frac{1}{|I_j(N_r)|}\sum_{n\in I_j(N_r)} x^{(j)}_n\Big)
		\ \xrightarrow[r\to\infty]{a.s.}\ \sum_{j=1}^{r}\lambda_j a_j.
	\end{equation}
	From the choice of the rational convex combination in \eqref{eq:xi_b_bounded} with \(\varepsilon>0\) is arbitrary and \eqref{eq:blocksum-subseq}, we conclude
	\begin{equation}\label{eq:blocksum-subseq-d}
	\frac{1}{N_r}\sum_{n=1}^{N_r}x_n
	\;\xrightarrow[r\to\infty]{\ \|\cdot\|_{\mathfrak X},\,a.s.\ }\; d .
	\end{equation}
	
	\medskip\noindent Now we extend \eqref{eq:blocksum-subseq-d} from the subsequence $N_r$ to all large $N$.
	Fix $a^\circ\in A$. For $N\in[N_r,N_r+q)$ define
	\[
	z_n:=
	\begin{cases}
		x_n, & 1\le n\le N_r,\\[2pt]
		x_n^{(a^\circ)}, & N_r<n\le N,
	\end{cases}
	\]
	where $x_n^{(a^\circ)}\in S_{X_n}(\mathcal F)$ are square–integrable selections with
	$E[x_n^{(a^\circ)}]=a^\circ$ as in (ii). Write $N=N_r+r$ with $0\le r<q$. Then
	\[
	\frac1N\sum_{n=1}^N z_n
	=\frac{N_r}{N}\Big(\frac1{N_r}\sum_{n=1}^{N_r}x_n\Big)
	+\frac1N\sum_{n=N_r+1}^{N}x_n^{(a^{\circ})} .
	\]
	For the remainder we estimate
	\[
	\Big\|\frac1N\sum_{n=N_r+1}^{N}x_n^{(a^\circ)}\Big\|_{\XX}
	\le \frac{r}{N}\max_{\,N_r<n\le N}\|x_n^{(a^\circ)}\|_{\XX}
	\le \frac{q}{N}\max_{\,N_r<n\le N}\|x_n^{(a^\circ)}\|_{\XX}.
	\]
	By (ii) and Borel–Cantelli again, $\|x_n^{(a^\circ)}\|_{\XX}/n\to0$ a.s., hence for any $\varepsilon>0$
	and all $N$ large enough,
	$\max_{N_r<n\le N}\|x_n^{(a^\circ)}\|_{\XX}\le \varepsilon N$, which yields
	\[
	\Big\|\frac1N\sum_{n=N_r+1}^{N}x_n^{(a^\circ)}\Big\|_{\XX}\le q\varepsilon\ \xrightarrow[N\to\infty]{}\ 0
	\quad \text{a.s.}
	\]
	Since $\frac{N_r}{N}\to1$ and $\frac1{N_r}\sum_{n=1}^{N_r}x_n\to d$ a.s. by \eqref{eq:blocksum-subseq-d}, we conclude
	$\frac1N\sum_{n=1}^N z_n\to d$ a.s., with $\frac1N\sum_{n=1}^N z_n\in S_N$.
	Therefore $d\in \mathrm{s}\text{-}\liminf_{n\to\infty}S_n$.
			
	Because $d\in D$ were arbitrary, we obtain
	\[
	D\subset\mathrm{s}\text{-}\liminf_{n\to\infty}S_n \quad\text{a.s.}
	\]

	Next, we prove that $\displaystyle\mathrm{w}\text{-}\limsup_{n\to\infty} S_n \subset D$ a.s. \\
	Let $x\in \displaystyle\mathrm{w}\text{-}\limsup_{n\to\infty} S_n$, so there exist $n_k\uparrow\infty$
	and $y_k\in S_{n_k}$ with $y_k\to x$ weakly. In particular $\{y_k\}$ is bounded.
	
	If $x\notin D$, by the strong Hahn–Banach separation for closed convex sets with
	recession we can choose $x^{*}\in S^{*}$ ($x^*\ne 0$) and $\varepsilon>0$ such that
	$s(x^{*},D)<\infty$ (so $x^*\in (V_{\infty}(D))^{\circ}$) and
	\[
	\langle x^{*},x\rangle \ \ge\ s(x^{*},D)+2\varepsilon .
	\]
	Since $y_k\to x$ weakly, $\langle x^{*},y_k\rangle\to\langle x^{*},x\rangle$, hence for $k$ large,
	\[
	\langle x^{*},y_k\rangle \ \ge\ s(x^{*},D)+\varepsilon . \tag{$\ast$}
	\]
	Using the basic identities of support functions for Minkowski sums and positive homogeneity (and the fact that $s(x^*,\operatorname{cl}B)=s(x^*,B)$), we have
	\[
	\begin{aligned}
		s(x^*,S_n)
		&= s\!\left(x^*, \frac1n\,\operatorname{cl}\sum_{k=1}^{n}X_k\right)
		\;=\; \frac1n\, s\!\left(x^*, \sum_{k=1}^{n}X_k\right)
		\;=\; \frac1n\sum_{k=1}^{n}s(x^*,X_k).
	\end{aligned}
	\]
	By weak stationarity and the Aumann–expectation identity,
	$E[s(x^*,X_n)]=s(x^*,E[X_n])=s(x^*,A)$ for all $n$ (cf. Remark~\ref{re:3.1}).
	Set
	\[
	Y_n:=s(x^*,X_n)-s(x^*,A),\qquad n\ge1 .
	\]
	Then $\{Y_n\}$ is a $\varphi$-mixing sequence (because $\sigma(Y_n)\subseteq\sigma(X_n)$),
	with $E[Y_n]=0$ and, by (iii),
	\(
	\sum_{n\ge1}\frac{E\,Y_n^{2}}{n^{2}}<\infty.
	\)
	Hence, applying Theorem~\ref{theo:1.1} with $a_n=n$ to the sequence $\{Y_n\}$ yields
	\[
	\frac1n\sum_{k=1}^{n}Y_k \;\longrightarrow\; 0 \quad\text{a.s.}
	\]
	Equivalently,
	\[
	s(x^*,S_n)=\frac1n\sum_{k=1}^{n}s(x^*,X_k)
	\;\longrightarrow\; s(x^*,A)\quad\text{a.s.}
	\]
	Finally, since $D=\overline{\operatorname{co}}A$ and the support function is
	unchanged by taking closed convex hull,
	\(
	s(x^*,D)=s(x^*,A)
	\)
	(possibly $+\infty$; in our case it is finite by assumption on $x^*$).
	Therefore $s(x^*,S_n)\to s(x^*,D)$ a.s.
	Therefore, for $k$ large, $s(x^{*},S_{n_k})\le s(x^{*},D)+\varepsilon/2$, which contradicts
	$(\ast)$ because $\langle x^{*},y_k\rangle\le s(x^{*},S_{n_k})$. Hence $x\in D$ and
	$\mathrm{w}\text{-}\limsup_{n\to\infty} S_n\subset D$ a.s.
	
\end{proof}

\begin{remark}[Functionals with infinite support value]\label{rem:inf-support}
	Let $A$ be (possibly) unbounded and set $D:=\overline{\operatorname{co}}A$. Recall that the
	support function $s(x^*,A)$ is finite exactly on the polar cone $C^\circ$ of
	the recession cone $C:=V_\infty(D)$, i.e.
	{ \[
	s(x^*,A)<\infty \quad\Longleftrightarrow\quad x^*\in C^\circ
	:=\{x^*\in \mathfrak{X}^*: \langle x^*,c\rangle\le 0\ \ \forall\,c\in C\}.
	\] }
	Hence assumption \textup{(iii)} of Theorem~\ref{theo:3.3} is only relevant
	for $x^*\in C^\circ$. On { $\mathfrak{X}^*\setminus C^\circ$ } we have $s(x^*,A)=+\infty$, so
	\textup{(iii)} is vacuous there.
	
	Two consequences:
	\begin{enumerate}\itemsep=0.2em
		\item If $s(x^*,A)=+\infty$ for all $x^*\in S^*$ (equivalently $C=\XX$ and thus
		$D=\XX$), then \textup{(iii)} is empty and Theorem~\ref{theo:3.3} holds under
		\textup{(i)}–\textup{(ii)} alone. In this case $S_n \xrightarrow{K\text{-}M} D=\XX$
		is trivial and the proof needs no change.
		\item In the proof of Theorem~\ref{theo:3.3} (the step $\mathrm{w}\text{-}\limsup_{n\to\infty} S_n\subset D$),
		the separating functional supplied by Hahn–Banach necessarily satisfies
		$s(x^*,D)<\infty$, hence $x^*\in C^\circ$. Therefore the use of \textup{(iii)}
		concerns only such $x^*$ and removing \textup{(iii)} outside $C^\circ$ does not
		affect the argument.
	\end{enumerate}
	
	\noindent\textit{Example (ray in $\R^2$).}
	For $A=\{(t,0):t\ge0\}$ one has $C=V_\infty(D)=\R_+(1,0)$ and {
	$C^\circ=\{x^*=(\xi_1,\xi_2)\in \mathfrak{X}^*:\ \xi_1\le 0\}$ }. Thus $s(x^*,A)=+\infty$
	for $\xi_1>0$, and \textup{(iii)} needs to be checked only for $\xi_1\le0$.
	This aligns with how the proof invokes \textup{(iii)}.
\end{remark}

The following example show that all assumptions of Theorem \ref{theo:3.3} do not force a degenerate situation.

\begin{example}[``Needle + shrinking halo'']\label{ex:3.5}
	Let $E=\mathbb R^2$ (Euclidean norm) and
	\[
	A:=\{(t,0):t\ge0\},\qquad D:=\overline{\operatorname{co}}\,A=A.
	\]
	Let $(\varepsilon_n)_{n\ge1}$ be i.i.d.\ uniform on $B(0,1)$ and put
	$z_n:=\varepsilon_n/n$. Define $X_n:=A\cup\{z_n\}$.
	
	\medskip\noindent\emph{Exact expansion of $\sum_{k=1}^n X_k$ and formula for $S_n$:}
	For $I\subseteq\{1,\dots,n\}$ write $Z_I:=\sum_{i\in I}z_i$ (with $Z_\varnothing:=0$).
	Using the distributive law of Minkowski sum over unions,
	$(U\cup V)+W=(U+W)\cup(V+W)$, together with $A+A=A$, we obtain
	\[
	\sum_{k=1}^n (A\cup\{z_k\})
	= \bigcup_{\substack{I\subseteq\{1,\dots,n\}\\ I\neq\{1,\dots,n\}}} \bigl(A+Z_I\bigr)
	\ \ \cup\ \ \{\,Z_{\{1,\dots,n\}}\,\}.
	\]
	Therefore
	\begin{equation}\label{eq:Sn-exact}
		S_n
		=\frac1n\,\operatorname{cl}\sum_{k=1}^n X_k
		= \Biggl(\ \bigcup_{\substack{I\subseteq\{1,\dots,n\}\\ I\neq\{1,\dots,n\}}}
		\Bigl(A+\tfrac1n Z_I\Bigr)\ \Biggr)\ \cup\ \Bigl\{\tfrac1n Z_{\{1,\dots,n\}}\Bigr\}.
	\end{equation}
	In particular, taking $I=\varnothing$ gives $A=\frac1n(A+Z_\varnothing)\subset S_n$, hence
	\[
	A\subset S_n\qquad\text{for all }n\ge1.
	\]
	
	\medskip\noindent\emph{A convenient upper inclusion and the ``halo'' bound:}
	From \eqref{eq:Sn-exact} and due to $\tfrac1n Z_{\{1,\dots,n\}}\subset A+\tfrac1nZ_{\{1,...,n\}}$ we have
	\[
	S_n\ \subset\bigcup_{I\subseteq\{1,...,n\}}\Bigl(A+\tfrac1n Z_I\Bigr)=A+\Bigl\{\tfrac1n Z_I:\ I\subseteq\{1,\dots,n\}\Bigr\}.
	\]
	Consequently,
	\[
	\max_{I\subseteq\{1,\dots,n\}}\Bigl\|\tfrac1n Z_I\Bigr\|
	\ \le\ \frac1n\sum_{i=1}^{n}\|z_i\|
	\ \le\ \frac1n\sum_{i=1}^{n}\frac1i
	=: r_n
	=O\!\Big(\frac{\log n}{n}\Big)\xrightarrow[n\to\infty]{}0.
	\]
	Thus $S_n\subset A+B(0,r_n)$ with $r_n\downarrow0$; geometrically, $S_n$ is $A$
	plus a shrinking bounded ``halo''.
	
	\medskip\noindent\emph{Recession cones:}
	Since $A\subset S_n$ we have $V_\infty(A)\subset V_\infty(S_n)$; and because
	$S_n\subset A+B(0,r_n)$ with $r_n<\infty$, it follows $V_\infty(S_n)\subset V_\infty(A)$.
	Hence
	\[
	V_\infty(S_n)=V_\infty(A)=\{t(1,0):t\ge0\}\qquad(\forall n).
	\]
	
	\medskip\noindent\emph{Verification of Theorem~\ref{theo:3.3}:}
	\begin{itemize}
		\item[(i)] $(X_n)$ are independent, so $\varphi(n)=0$ and $\sum_n\varphi(n)^{1/2}<\infty$.
		\item[(ii)] Fix $a=(t_0,0)\in A$. Because $A\subset X_n(\omega)$ for every $n$ and $\omega$,
		the constant selection $x_n(\omega)\equiv a$ lies in $S_{X_n}^2$,
		satisfies $E[x_n]=a$, and
		\[
		\sum_{n=1}^{\infty}\frac{E\|x_n-a\|^2}{n^2}=0<\infty.
		\]
		(When $a=0$, one may also take $x_n=z_n$, with $E\|x_n\|^2\le n^{-2}$ so that
		$\sum n^{-4}<\infty$.)
		\item[(iii)] For any $x^*=(x_1^*,x_2^*)\in S^*$ with $s(x^*,A)<\infty$ we have $x_1^*\le0$ and
		$s(x^*,A)=0$. Then
		\[
		0\le s(x^*,X_n)-s(x^*,A)
		=\max\{\langle x^*,z_n\rangle,0\}
		\le |\langle x^*,z_n\rangle|
		\le \|z_n\|.
		\]
		Hence $E|s(x^*,X_n)-s(x^*,A)|^2\le E\|z_n\|^2\le n^{-2}$ and
		\[
		\sum_{n=1}^{\infty}\frac{E|s(x^*,X_n)-s(x^*,A)|^2}{n^2}
		\ \le\ \sum_{n=1}^{\infty}\frac{1}{n^4}\ <\ \infty,
		\]
		which is precisely the pointwise summability required by (iii).
	\end{itemize}
	
	\medskip\noindent\emph{Kuratowski--Mosco convergence:}
	Because $S_n\subset A+B(0,r_n)$ with $r_n\to0$ and $A\subset S_n$ for all $n$,
	we have $\mathrm{w}\text{-}\limsup_{n\to\infty}S_n
	=\mathrm{s}\text{-}\liminf_{n\to\infty}S_n=A=D$.
	Thus $S_n\xrightarrow{\mathrm{K\text{-}M}}D$.
\end{example}
 Our next example shows that if the pointwise summability in \textup{(iii)} is dropped
 (while retaining \textup{(i)}--\textup{(ii)}), the conclusion of Theorem~\ref{theo:3.3}
 can fail. In that sense, \textup{(iii)} is genuinely needed in the unbounded case.

\begin{example}\label{ex:3.6}
	(\emph{Violating (iii): the K--M limit may fail.})
	Let $E=\mathbb R^2$ with the Euclidean norm and let
	\[
	A:=\{(t,0):t\ge 0\},\qquad D:=\overline{\operatorname{co}}\,A=A.
	\]
	For each $n\ge1$ let $\theta_n$ take the values $\pm \tfrac1n$ with probability
	$1/2$ each, independently across $n$, and set
	\[
	v_n:=(\cos\theta_n,\ \sin\theta_n),\qquad
	X_n:=\{\,t\,v_n:\ t\ge 0\,\}\quad(\text{a closed ray from the origin}).
	\]
	Thus each $X_n$ is a closed convex cone (a ray). Define
	\[
	S_n:=\frac1n\,\operatorname{cl}\sum_{k=1}^n X_k,\qquad
	V_\infty(C):=\{u\in E:\ C+\lambda u\subset C\ \text{for all }\lambda\ge 0\}.
	\]
	
	\smallskip\noindent\emph{(i) and (ii) hold:}
	The sequence $(X_n)$ is independent, hence $\varphi(n)=0$ for all $n$
	and $\sum_n \varphi(n)^{1/2}<\infty$.
	Fix $a=(a,0)\in A$ and set $t_n:=a/\cos(1/n)$.
	Define a selection $x_n(\omega):=t_n\,v_n(\omega)\in X_n(\omega)$.
	Then
	\[
	E[x_n]=\tfrac12 t_n(\cos\tfrac1n,\sin\tfrac1n)+
	\tfrac12 t_n(\cos\tfrac1n,-\sin\tfrac1n)
	=(a,0)=a,
	\]
	and
	\[
	\|x_n-a\|^2
	= a^2\tan^2\!\tfrac1n
	\sim \frac{a^2}{n^2}\quad(n\to\infty).
	\]
	Hence
	\[
	\sum_{n=1}^\infty \frac{E\|x_n-a\|^2}{n^2}
	\ \asymp\ \sum_{n=1}^\infty \frac{a^2}{n^4}\ <\ \infty,
	\]
	so (ii) is satisfied.
	
	\smallskip\noindent\emph{(iii) fails:}
	Take $x^*:=(0,1)\in E^*$. Then $s(x^*,A)=0$, while
	\[
	s(x^*,X_n)=\sup_{t\ge 0} t\,\langle x^*,v_n\rangle
	=\begin{cases}
		+\infty,& \theta_n=+1/n,\\
		0,& \theta_n=-1/n,
	\end{cases}
	\]
	so $E|\,s(x^*,X_n)-s(x^*,A)\,|^2=+\infty$ for every $n$ and
	$\sum_n \frac{E|\,s(x^*,X_n)-s(x^*,A)\,|^2}{n^2}=+\infty$.
	Thus assumption \textup{(iii)} is violated.
	
	\smallskip\noindent\emph{No Kuratowski--Mosco convergence to $D$:}
	Each $X_k$ is the cone $\mathrm{cone}\{v_k\}$; the Minkowski sum of finitely
	many cones is the cone generated by the union of their generators, hence
	\[
	\sum_{k=1}^n X_k
	=\mathrm{cone}\{v_1,\dots,v_n\},
	\qquad
	S_n=\tfrac1n\,\mathrm{cone}\{v_1,\dots,v_n\}
	=\mathrm{cone}\{v_1,\dots,v_n\}.
	\]
	Let
	\[
	\alpha_n^+:=\max\{\theta_k:\ 1\le k\le n,\ \theta_k>0\},\qquad
	\alpha_n^-:=\min\{\theta_k:\ 1\le k\le n,\ \theta_k<0\}.
	\]
	Then $S_n$ is the closed convex cone bounded by the two rays of directions
	$\alpha_n^+$ and $\alpha_n^-$:
	\[S_n=\mathrm{cone}\big\{(\cos\alpha_n^+,\sin\alpha_n^+),\ (\cos\alpha_n^-,\sin\alpha_n^-)\big\}.\]
	
	In particular $V_\infty(S_n)=\mathrm{cone}\{v_k:\ 1\le k\le n\}$.
	By independence, with probability one there exist finite indices
	$k_+\!,k_->0$ such that $\theta_{k_+}=+1/k_+$ and $\theta_{k_-}=-1/k_-$.
	Hence, for all $n\ge N_0:=\max\{k_+,k_-\}$,
	\[
	S_n
	\supset \mathrm{cone}\big\{(\cos \tfrac1{k_+},\sin \tfrac1{k_+}),
	(\cos \tfrac1{k_-},-\sin \tfrac1{k_-})\big\},
	\]
	which is a fixed sector of positive opening angle, strictly larger than $D=A$.
	Therefore $\mathrm{s}\text{-}\liminf_{n\to\infty}S_n$ contains that sector and
	$\mathrm{w}\text{-}\limsup_{n\to\infty}S_n$ does as well; in particular
	\[
	\mathrm{w}\text{-}\limsup_{n\to\infty}S_n\ \not\subset\ D,
	\qquad
	S_n\ \text{does not converge to }D\text{ in the K--M sense}.
	\]
\end{example}

\section{Conclusions}\label{conc}

In this paper, we have extended the definitions of weak stationarity and
$\varphi$-mixing to set-valued random variables and, under natural
summability hypotheses, established several strong laws of large
numbers in both the Hausdorff and Kuratowski-Mosco frameworks.
The illustrative examples show that our conditions neither force
degeneracy to single-valued variables nor are they logically
necessary-only sufficient.

Future work may include
\emph{(i)} extending the results to stronger mixing concepts such as
$\rho$-mixing, whose dependence coefficient is defined by
\[
\rho(\mathcal A,\mathcal B):=
\sup_{X\in\mathcal U(\mathcal A),\,Y\in\mathcal U(\mathcal B)}
\Bigl|\frac{\operatorname{Cov}(X,Y)}
{\sqrt{\operatorname{Var}(X)\operatorname{Var}(Y)}}\Bigr|
\quad\text{\cite{yang2005dp}};
\]
and \emph{(ii)} proving strong laws for weakly stationary
$\varphi$-mixing set-valued sequences by means of summability methods as
in~\cite{kiesel1997strong}.

\section*{Acknowledgements}
We thank the anonymous referees for their valuable comments and insightful suggestions, which greatly improved the manuscript.

\section*{Funding}
This work was supported by the Institute of Information Technology,
Vietnam Academy of Science and Technology, under project
CSCL02.03/25-26.

\section*{ORCID}
Luc Tri Tuyen: 0000-0002-4822-2978

\section*{Declaration of generative AI and AI-assisted technologies in the writing process}

During the preparation of this work the author used ChatGPT in order to improve language and readability. After using this tool/service, the author reviewed and edited the content as needed and take full responsibility for the content of the publication.

\bibliographystyle{plain}      
\bibliography{references_3}   
\end{document}